\documentclass[11pt]{article}
\usepackage{amsfonts,latexsym,amsmath,amscd,geometry}
\usepackage{amssymb}
\usepackage{latexsym}

\newcommand \nc{\newcommand}
\newtheorem{theorem}{Theorem}[section]
\newtheorem{lemma}[theorem]{Lemma}

\newtheorem{corollary}[theorem]{Corollary}

\nc{\ba}{\begin{array}}\nc{\ea}{\end{array}}
\nc{\be}{\begin{eqnarray}}\nc{\ee}{\end{eqnarray}}
\nc{\beq}{\begin{equation}}\nc{\eeq}{\end{equation}}
\nc{\bex}{\begin{eqnarray*}}\nc{\eex}{\end{eqnarray*}}
\nc{\btm}{\begin{theorem}} \nc{\etm}{\end{theorem}}
\nc{\blm}{\begin{lemma}} \nc{\elm}{\end{lemma}}
\nc{\R}{\mathbb{R}} \nc{\va}{\varepsilon} \nc{\ls}{\limits}

\def\o{\omega}\def\tr{\triangle}
\def\pf{\noindent{\bf Proof.\quad}}\def\endpf{\hfill$\Box$}
\nc \qed {\hfill $\Box$}

\begin{document}
\title{Blow up criterion for nematic liquid crystal flows}
\author{Tao Huang\footnote{Department of Mathematics, University of Kentucky,
Lexington, KY 40506. The work is partially supported by NSF 1000115.} \quad Changyou Wang$^*$ }
\date{}
\maketitle

\begin{abstract}
 In this paper, we establish a blow up criterion for the short time classical solution of the nematic liquid crystal flow, a simplified version of Ericksen-Leslie system modeling the hydrodynamic evolution of nematic liquid crystals, in dimensions two and three. More precisely, $0<T_*<+\infty$ is the maximal time interval iff
(i) for $n=3$, $|\o|+|\nabla d|^2\notin L^{1}_tL^{\infty}_x(\mathbb R^3\times [0,T_*])$;
and (ii) for $n=2$, $|\nabla d|^2\notin L^1_t L^\infty_x(\mathbb R^2\times [0,T_*])$.
\end{abstract}

%\vskip 5mm
\section {Introduction}
\setcounter{equation}{0}
\setcounter{theorem}{0}

In this paper, we consider the Cauchy problem to the following hydrodynamic flow of nematic liquid crystals in $\R^n$ ($n=2$ or $3$):
\begin{align}
u_t+u\cdot\nabla u-\nu\triangle u+\nabla p&=-\Delta d\cdot \nabla d,\label{liquidcrystal-1}\\
\nabla \cdot u&=0, \label{liquidcrystal-2} \\
d_t+u\cdot\nabla d&=\triangle d+|\nabla d|^2d,\label{liquidcrystal-3}\\
(u, d)\big|_{t=0}&= (u_0, d_0), \label{liquidcrystal-4}
\end{align}
where $u:\R^n\rightarrow\R^n$ represents the velocity field of the incompressible viscous fluid, $\nu>0$ is the Kinematic viscosity, $p:\R^3\rightarrow\R$ represents the pressure function, $d:\R^n\rightarrow \mathbb{S}^2$ represents the macroscopic average of the nematic liquid crystal orientation field, $\nabla\cdot$ and $\Delta$ denotes the divergence operator and the Laplace operator respectively,
$u_0:\mathbb R^n\to\mathbb R^n$ is a given initial velocity field with $\nabla\cdot u_0=0$,
and $d_0:\mathbb R^n\to S^2$ is a given initial liquid crystal orientation field.

The  system (\ref{liquidcrystal-1})-(\ref{liquidcrystal-3}) is a simplified version of the Ericksen-Leslie system modeling the hydrodynamics of nematic liquid crystals developed during the period of 1958 through 1968 (\cite{Gennes} \cite{Ericksen} \cite{Leslie}). It is
a macroscopic continuum description of the time evolution of the material under
the influence of both the flow field $u(x,t)$, and the macroscopic description of the
microscopic orientation configurations $d(x,t)$ of rod-like liquid crystals. Recall that the Ericksen-Leslie theory  reduces to the Ossen-Frank theory in the static case, see Hardt-Lin-Kinderlehrer \cite{HLK} and references therein.
The system (\ref{liquidcrystal-1})-(\ref{liquidcrystal-3}) was first introduced by Lin and Liu in their important works \cite{LL, LL1} during the 1990's. Roughly speaking, (\ref{liquidcrystal-1})-(\ref{liquidcrystal-3}) is a system that couples between the non-homogeneous Navier-Stokes equation and the transported heat flow of harmonic maps into $S^2$. For dimension $n=2$,  Lin-Lin-Wang \cite{Lin-Lin-Wang}
have proved the global existence of Leray-Hopf type weak solutions to (\ref{liquidcrystal-1})-(\ref{liquidcrystal-3}) on bounded domains in $\mathbb R^2$ under the initial and boundary value conditions (see \cite{Hong} for the case $\Omega=\mathbb R^2$), and Lin-Wang \cite{LW} have further established the uniqueness for such weak solutions. It is an interesting and challenging problem to study the nematic liquid crystal flow equation (\ref{liquidcrystal-1})-(\ref{liquidcrystal-3}) in dimension three,
such as the global existence of weak solutions and the partial regularity of suitable weak solutions.

In this paper, we will consider the short time classical solution to (\ref{liquidcrystal-1}) -(\ref{liquidcrystal-4})
and  address some criterion that characterizes the first finite singular time. It is well-known that if the initial velocity $u_0\in H^s(\mathbb R^n,
\mathbb R^n)$ with $\nabla\cdot u_0=0$ and $d_0\in H^{s+1}(\mathbb R^n,S^2)$ for $s\geq n$, then there is  $T_0>0$ depending only on $\|u_0\|_{H^s}$ and $\|d_0\|_{H^{s+1}}$ such that (\ref{liquidcrystal-1})-(\ref{liquidcrystal-4}) has a unique, classical solution $(u,d)$ in $\R^n\times[0,T_0)$ satisfying
\begin{align}
\label{blp1.2} & u\in C([0,T],H^s(\mathbb R^n))\cap C^1([0,T], H^{s-1}(\mathbb R^n)) \mbox{ and }\nonumber\\
&
d\in C([0,T],H^{s+1}(\mathbb R^n,S^2))\cap C^1([0,T], H^{s}(\mathbb R^n,S^2)),
\end{align}
for any $0<T<T_0$.  Assume $T_*>0$ is the maximum value such that (1.5) holds with $T_0=T_*$. We would like to characterize such a
$T_*$.

Recall that when $d$ is a constant vector, (\ref{liquidcrystal-1})-(\ref{liquidcrystal-4}) becomes the Navier-Stokes equation.
In their famous work \cite{Beale-Kato-Majda}, Beale-Kato-Majda proved that for $n=3$, if $T_*>0$ is the first finite singular time,
then the vorticity $\omega=\nabla\times u$ doesn't  belong to $L^1_tL^\infty_x(\mathbb R^3\times [0,T_*))$.
On the other hand, when $u=0$, (\ref{liquidcrystal-1})-(\ref{liquidcrystal-4}) becomes the heat flow of harmonic maps into
$S^2$,  Wang proved in \cite{Wang} that for $n\ge 2$, if $\nabla d\in L^\infty_t L^n_x(\mathbb R^n\times [0,T])$, then
$d\in C^\infty(\mathbb R^n\times (0,T])$.  Our main result on (\ref{liquidcrystal-1})-(\ref{liquidcrystal-4}) is a natural extension
of \cite{Beale-Kato-Majda} and \cite{Wang}.
\btm\label{maintheorem}{\it For $n=3$, $s\ge 3$, $u_0\in H^s(\mathbb R^n,\mathbb R^n)$ with $\nabla\cdot u_0=0$ and
$d_0\in H^{s+1}(\mathbb R^n,S^2)$, let $T_*>0$ be the maximum value such that (\ref{liquidcrystal-1})-(\ref{liquidcrystal-4})
has a unique solution $(u,d)$ satisfying (1.5) with $T_0$ replaced by $T_*$. If $T_*<+\infty$, then
\beq\label{blpcondition}
\int_{0}^{T_*}\left(\|\o(t)\|_{L^{\infty}(\R^3)}+\|\nabla d(t)\|^2_{L^{\infty}(\R^3)} \right)\,dt=\infty,
\eeq
where $\o=\nabla\times u$ is the vorticity. In particular,
\beq\label{blp1.3}
\limsup\ls_{t\nearrow T_*}\left(\|\o(t)\|_{L^{\infty}(\R^3)}+\|\nabla d(t) \|^2_{L^{\infty}(\R^3)}\right)=\infty.
\eeq
}\etm

As a byproduct of the proof of theorem 1.1 and the regularity theorem by \cite{Lin-Lin-Wang}, we obtain a corresponding criterion in dimension $n=2$.
More precisely, we have

\begin{corollary}\label{2-d case}{\it  For $n=2$, $s\ge 2$, $u_0\in H^s(\mathbb R^n,\mathbb R^n)$ with $\nabla\cdot u_0=0$ and
$d_0\in H^{s+1}(\mathbb R^n,S^2)$, let $T_*>0$ be the maximum value such that (\ref{liquidcrystal-1})-(\ref{liquidcrystal-4})
has a unique solution $(u,d)$ satisfying (1.5) with $T_0$ replaced by $T_*$. If $T_*<+\infty$, then
\beq\label{blpcondition2-d}
\int_{0}^{T_*}\|\nabla d(t)\|^2_{L^{\infty}(\R^2)} \,dt=\infty.
\eeq
In particular,
\beq\label{blp1.5}
\limsup\ls_{t\nearrow T_*}\|\nabla d(t) \|_{L^{\infty}(\R^2)}=\infty.
\eeq
}\end{corollary}

\section {Proof of Theorem 1.1}
\setcounter{equation}{0}
\setcounter{theorem}{0}

For simplicity, we assume $\nu=1$. We need the following lemma to prove theorem 1.1.

\blm\label{blplemma2.1}{\it For $n=2$ or $3$, $s\ge n$, $u_0\in H^s(\mathbb R^n,\mathbb R^n)$ with $\nabla\cdot u_0=0$ and
$d_0\in H^{s+1}(\mathbb R^n,S^2)$, $M>0$, and $T_0>0$,
let $(u,d)$ be a solution to (\ref{liquidcrystal-1})-(\ref{liquidcrystal-4}) satisfying (\ref{blp1.2}) and
\begin{equation}\label{blp2.1}
\begin{cases}\int_{0}^{T_0}\left(\|\o(t)\|_{L^{\infty}(\R^n)}+\|\nabla d(t)\|^2_{L^{\infty}(\R^n)} \right)\,dt\leq M & {\rm{ for }}\ n=3,\\
{\rm{or}}\ \ \ \ \ \ \ \ \ \ \ \ \ \ \ \ \ \ \ \ \int_{0}^{T_0}\|\nabla d(t)\|^2_{L^{\infty}(\R^n)} \,dt\leq M&  {\rm{ for }}\ n=2.\end{cases}
\end{equation}
Then
\beq\label{blp2.2}
\sup\ls_{0\leq t\leq T_0}\left(\|\o(t)\|_{L^{2}(\R^n)}+\|\nabla^2 d(t)\|_{L^{2}(\R^n)} \right)\leq C,
\eeq
where $C>0$ depends only on $u_0$, $d_0$ and $M$.
}\elm

\pf  Taking $\nabla\times$ on (\ref{liquidcrystal-1}), we obtain
\begin{equation}\label{blp2.3}
\o_t-\triangle\o+u\cdot\nabla\o=\begin{cases}
\o\cdot\nabla u-\nabla\times(\triangle d\cdot\nabla d) &\ {\rm{if}}\ n=3,\\
-\nabla\times(\triangle d\cdot\nabla d) &\ {\rm{if}}\ n=2
\end{cases}
\end{equation}
Multiplying $\o$ and integrating over $\R^n$, we obtain
\beq\label{blp2.4}
\frac{1}{2}\frac{d}{dt}\int_{\R^n}|\o|^2\,dx+\int_{\R^n}|\nabla \o|^2\,dx=
\begin{cases}\int_{\R^n}[(\o\cdot\nabla)u\cdot\o +(\triangle d\cdot\nabla d)\cdot(\nabla\times)\o]\,dx
& {\rm{for}}\ n=3,\\
\int_{\R^n}(\triangle d\cdot\nabla d)\cdot(\nabla\times\o) \,dx & {\rm{for}}\ n=2.
\end{cases}
\eeq
where we have used the fact
$$\int_{\R^n}(u\cdot\nabla)\o\cdot\o \,dx=\frac12\int_{\R^n}(u\cdot\nabla)|\o|^2 \,dx=0.$$
Since
$$\nabla u=(-\triangle)^{-1}\nabla(\nabla\times\o),$$
we have $\|\nabla u\|_{L^2}\leq C\|\o\|_{L^2}$ and
\beq\label{blp2.5}
\left|\int_{\R^n}(\o\cdot\nabla)u\cdot\o \,dx\right|\leq C\|\o\|_{L^{\infty}}\|\o\|_{L^2}^2.
\eeq
By Young's inequality, we obtain
\beq\label{blp2.6}
\left|\int_{\R^n}(\triangle d\cdot\nabla d)\cdot(\nabla\times\o) \,dx\right|\leq C\int_{\R^n}|\triangle d|^2|\nabla d|^2\,dx+\frac{1}{2}\int_{\R^n}|\nabla\o|^2\,dx.
\eeq
Combining (\ref{blp2.4}), (\ref{blp2.5}), with (\ref{blp2.6}) , we have
\beq\label{blp2.7}
\frac{d}{dt}\|\o\|_{L^2}^2+\|\nabla \o\|_{L^2}^2\leq
\begin{cases} C\|\o\|_{L^{\infty}}\|\o\|_{L^2}^2+C\int_{\R^n}|\triangle d|^2|\nabla d|^2\,dx & {\rm{for}}\ n=3,\\
C\int_{\R^n}|\triangle d|^2|\nabla d|^2\,dx & {\rm{for}}\ n=2.
\end{cases}
\eeq

Taking $\triangle$ on (\ref{liquidcrystal-3}), multiplying $\triangle d$ and integrating over $\R^n$, we obtain
\beq\label{blp2.8}
\frac{1}{2}\frac{d}{dt}\|\tr d\|_{L^2}^2+\|\nabla\tr d\|_{L^2}^2= \int_{\R^n}\triangle(|\nabla d|^2d)\cdot\tr d \,dx
-\int_{\R^n}\triangle(u\cdot\nabla d)\cdot\tr d \,dx.
\eeq
Since
$$\int_{\R^n}(u\cdot\nabla\tr d)\cdot\tr d\,dx=\frac12\int_{\R^n}(u\cdot\nabla)(|\tr d|^2)\,dx=0,$$
and
$$\nabla\times\o=\nabla\times(\nabla\times u)=\nabla(\nabla\cdot u)-\tr u=-\tr u,$$
we obtain,
\beq\label{blp2.9}\begin{split}
\left|\int_{\R^n}\triangle(u\cdot\nabla d)\cdot\tr d \,dx\right|&\leq\int_{\R^n}|\tr u||\nabla d||\tr d| \,dx
+2\int_{\R^n}|\nabla u||\nabla^2 d||\tr d| \,dx\\
&\leq\int_{\R^n}|\nabla \o||\nabla d||\tr d| \,dx+2\int_{\R^n}|\nabla u||\nabla^2 d||\tr d| \,dx\\
&=I_1+I_2.
\end{split}\eeq
\beq\label{blp2.10}
I_1\leq \frac{1}{4}\|\nabla \o\|_{L^2}^2+C\|\nabla d\|_{L^{\infty}}^2\|\tr d\|_{L^2}^2.
\eeq
\beq\label{blp2.11}
\begin{split}
I_2&\leq C \|\nabla u\|_{L^2}\|\nabla^2d\|_{L^4}^2\\
&\leq C\|\o\|_{L^2}\|\nabla d\|_{L^{\infty}}\|\nabla^3 d\|_{L^2}\\
&\leq C\|\o\|_{L^2}\|\nabla d\|_{L^{\infty}}\|\nabla\tr d\|_{L^2}\\
&\leq \frac{1}{4}\|\nabla\tr d\|_{L^2}^2+C\|\o\|_{L^2}^2\|\nabla d\|_{L^{\infty}}^2,
\end{split}\eeq
where we have used Nirenberg's interpolation inequality: for nonegative integers $k$ and $l$ with $k\leq l-1$,
$$\|\nabla^k f\|_{L^{\frac{2l}{k}}}^{\frac{2l}{k}}\leq C\|f\|_{L^{\infty}}^{\frac{2l}{k}-2}\|\nabla^l f\|_{L^2}^2.$$
Combining (\ref{blp2.9}), (\ref{blp2.10}), with (\ref{blp2.11}), we have
\beq\label{blp2.12}
\left|\int_{\R^n}\triangle(u\cdot\nabla d)\cdot\tr d \,dx\right|
\leq \frac{1}{4}\|\nabla \o\|_{L^2}^2+\frac{1}{4}\|\nabla\tr d\|_{L^2}^2+C\|\nabla d\|_{L^{\infty}}^2(\|\tr d\|_{L^2}^2+\|\o\|_{L^2}^2).
\eeq

Now we need to estimate the first term in the right hand side of  (\ref{blp2.8}).
\beq\label{blp2.13}
\begin{split}
&\int_{\R^n}\triangle(|\nabla d|^2d)\cdot\tr d \,dx\\
=&\int_{\R^n}\triangle(|\nabla d|^2)d\cdot\tr d \,dx+\int_{\R^n}|\nabla d|^2|\tr d|^2 \,dx
+\int_{\R^n}2\nabla|\nabla d|^2\cdot\nabla d\tr d \,dx\\
=&I_3+I_4+I_5.
\end{split}
\eeq
By integration by parts, we obtain
\beq\label{blp2.14}
\begin{split}
|I_3|&=\left|-\int_{\R^n}(\nabla(|\nabla d|^2)\nabla d\cdot\tr d+\nabla(|\nabla d|^2)d\cdot\nabla\tr d) \,dx\right|\\
&\leq 2\int_{\R^n}(|\nabla d|^2|\nabla^2 d||\tr d|+|\nabla d||\nabla^2d||\nabla\tr d|) \,dx\\
&\leq C\|\nabla d\|_{L^{\infty}}^2\|\tr d\|_{L^2}^2+\frac{1}{4}\|\nabla \tr d\|_{L^2}^2.
\end{split}\eeq
\beq\label{blp2.15}
|I_4|\leq \|\nabla d\|_{L^{\infty}}^2\|\tr d\|_{L^2}^2.\\
\eeq
\beq\label{blp2.16}
\begin{split}
|I_5|&\leq 4\int_{\R^n}|\nabla d|^2|\nabla^2 d||\tr d|\,dx\\
&\leq C\|\nabla d\|^2_{L^{\infty}}\|\nabla^2d\|_{L^2}\|\tr d\|_{L^2}\\
&\leq C\|\nabla d\|^2_{L^{\infty}}\|\tr d\|_{L^2}^2.
\end{split}\eeq
Combining (\ref{blp2.13}), (\ref{blp2.14}), (\ref{blp2.15}), with (\ref{blp2.16}), we have
\beq\label{blp2.17}
\int_{\R^3}\triangle(|\nabla d|^2d)\cdot\tr d \,dx
\leq C\|\nabla d\|_{L^{\infty}}^2\|\tr d\|_{L^2}^2+\frac{1}{4}\|\nabla \tr d\|_{L^2}^2.
\eeq
Combining (\ref{blp2.8}), (\ref{blp2.12}) and (\ref{blp2.17}), we obtain
\beq\label{blp2.18}
\frac{d}{dt}\|\tr d\|_{L^2}^2+\|\nabla\tr d\|_{L^2}^2
\leq  \frac{1}{2}\|\nabla \o\|_{L^2}^2+C\|\nabla d\|_{L^{\infty}}^2(\|\tr d\|_{L^2}^2+\|\o\|_{L^2}^2).
\eeq

Adding (\ref{blp2.7}) and (\ref{blp2.18}) together, we obtain
\beq\label{blp2.19}
\begin{split}
&\frac{d}{dt}(\|\o\|_{L^2}^2+\|\tr d\|_{L^2}^2)+\frac{1}{2}\|\nabla \o\|_{L^2}^2+\|\nabla\tr d\|_{L^2}^2\\
\leq& \begin{cases} C\left(\|\o\|_{L^{\infty}}+\|\nabla d\|_{L^{\infty}}^2\right)\left(\|\tr d\|_{L^2}^2+\|\o\|_{L^2}^2\right)
& {\rm{for}} \ n=3\\
C\|\nabla d\|_{L^{\infty}}^2\left(\|\tr d\|_{L^2}^2+\|\o\|_{L^2}^2\right)
& {\rm{for}}\ n=2
\end{cases}
\end{split}\eeq
Then by Gronwall's inequality,
\beq\label{blp2.20}
\begin{split}
&\|\o(T_0)\|_{L^2}^2+\|\tr d(T_0)\|_{L^2}^2\\
\leq &\begin{cases}\left(\|\tr d_0\|_{L^2}^2+\|\o_0\|_{L^2}^2\right)
\exp\left(C\int_0^{T_0}\left(\|\o(t)\|_{L^{\infty}}+\|\nabla d(t)\|_{L^{\infty}}^2\right)\,dt\right) & {\rm{for}}\ n=3\\
\left(\|\tr d_0\|_{L^2}^2+\|\o_0\|_{L^2}^2\right)
\exp\left(C\int_0^{T_0}\|\nabla d(t)\|_{L^{\infty}}^2\,dt\right) & {\rm{for}}\ n=2.
\end{cases}
\end{split}
\eeq
Since
$$\int_{\mathbb R^n}|\Delta d|^2\,dx=\int_{\mathbb R^n}|\nabla^2 d|^2\,dx,$$
this yields the conclusion and hence completes the proof of lemma \ref{blplemma2.1}.
\endpf

\vspace{5mm}
{\bf Proof of Theorem \ref{maintheorem}}: We prove the theorem by contradiction. Assume that (\ref{blpcondition}) were not true. Then
there is  $0<M<\infty$ such that
\beq\label{blp2.21}
\int_{0}^{T_*}\left(\|\o(t)\|_{L^{\infty}(\R^3)}+\|\nabla d(t)\|^2_{L^{\infty}(\R^3)} \right)dt\leq M.
\eeq
Then by lemma \ref{blplemma2.1}, we have
\beq\label{blp2.22}
\sup\ls_{0\leq t\leq T_*}\left(\|\o(t)\|_{L^{2}(\R^3)}+\|\nabla^2 d(t)\|_{L^{2}(\R^3)} \right)<C,
\eeq
where $C>0$ depends on $u_0$, $d_0$ and $M$.

 If we could control $\|u(t)\|_{H^3}+\|\nabla d(t)\|_{H^3}$ for any $0\leq t\leq T_*$ in  terms of $u_0$, $d_0$ and $M$,
we would reach a contradiction. To do this, we need higher order energy estimates, which can be done as follows.

For any multi-index $s$ with $|s|=3$, taking $D^s$ on (\ref{liquidcrystal-1}), multiplying $D^s u$ and integrating over $\R^3$, we obtain
\beq\label{blp2.23}
\begin{split}
&\frac{d}{dt}\int_{\R^3}\frac{|D^s u|^2}{2}\,dx+\int_{\R^3}|D^{s+1}u|^2\,dx\\
=&-\int_{R^3}\left[D^s(u\cdot\nabla u)-u\cdot\nabla D^su\right]\cdot D^su \,dx\\
&-\int_{\R^3} D^s(\tr d\cdot\nabla d)\cdot D^su \,dx\\
=&J_1+J_2.
\end{split}
\eeq
For $J_1$, we need to use the following estimate (see \cite{Kato-Ponce})
$$
\|D^s(fg)-fD^sg\|_{L^2}\leq C(\|f\|_{H^3}\|g\|_{L^{\infty}}+\|\nabla f\|_{L^{\infty}}\|g\|_{H^{2}}).
$$
Setting $f=u$ and $g=\nabla u$, we have
\beq\label{blp2.24}
\begin{split}
|J_1|
&\leq C\|D^s(u\cdot\nabla u)-u\cdot\nabla D^su\|_{L^2}\|D^su\|_{L^{2}}^2\\
&\leq C\|\nabla u\|_{L^{\infty}}\|u\|_{H^{3}}^2.
\end{split}
\eeq
Applying the Leibniz's rule and Nirenberg's interpolation inequality, we have
\beq\label{blp2.25}\begin{split}
J_2&=\int_{\R^3} D^{s-1}(\tr d\cdot\nabla d)\cdot D^{s+1}u \,dx\\
&\leq \frac{1}{2}\|D^{s+1}u\|_{L^2}^2+C\int_{\R^3}|D^{s-1}(\tr d\cdot\nabla d)|^2\,dx\\
&\leq \frac{1}{2}\|D^{s+1}u\|_{L^2}^2+C\int_{\R^3}(|\nabla^{4}d|^2|\nabla d|^2+|\nabla^{2}d|^2|\nabla^3 d|^2)\,dx\\
&\leq \frac{1}{2}\|D^{s+1}u\|_{L^2}^2+C\|\nabla d\|_{L^{\infty}}^2\|\nabla d\|_{H^3}^2
+C\|\nabla^{2}d\|_{L^{6}}^2\|\nabla^{3}d\|_{L^3}^2\\
&\leq \frac{1}{2}\|D^{s+1}u\|_{L^2}^2+C\|\nabla d\|_{L^{\infty}}^2\|\nabla d\|_{H^3}^2.
\end{split}
\eeq
Combining (\ref{blp2.23}), (\ref{blp2.24}), with (\ref{blp2.25}), we have
\beq\label{blp2.26}
\frac{d}{dt}\|D^s u\|_{L^2}^2+\|D^{s+1} u\|_{L^2}^2\leq C\left(\|\nabla u\|_{L^{\infty}}\|u\|_{H^3}^2+\|\nabla d\|_{L^{\infty}}^2\|\nabla d\|_{H^3}^2\right).
\eeq

Taking $D^{s+1}$ on (\ref{liquidcrystal-3}), multiplying $D^{s+1} d$ and integrating over $\R^3$, we obtain
\beq\label{blp2.27}
\begin{split}
&\frac{d}{dt}\int_{\R^3}\frac{|D^{s+1} d|^2}{2}\,dx+\int_{\R^3}|D^{s+2} d|^2\,dx\\
=&-\int_{R^3}\left(D^{s+1}(u\cdot\nabla d)-u\cdot\nabla D^{s+1}d\right)D^{s+1}d \,dx\\
&+\int_{\R^3} D^{s+1}(|\nabla d|^2d)\cdot D^{s+1}d \,dx\\
=&J_3+J_4.
\end{split}
\eeq
For $J_3$, similar as the proof of (\ref{blp2.24}), we have
\beq\label{blp2.28}
\begin{split}
J_3&\leq C\|D^{s+1}(u\cdot\nabla d)-u\cdot\nabla
D^{s+1}d\|_{L^2}\|D^{s+1}d\|_{L^2 }\\
&\leq C\|\nabla d\|_{L^{\infty}}\|u\|_{H^4}\|D^{s+1}d\|_{L^2
}+\|\nabla u\|_{L^{\infty}}\|\nabla d\|_{H^3}\|D^{s+1}d\|_{L^2 }\\
&\leq \frac{\epsilon}{2}\|u\|_{H^4}^2+C\|\nabla
d\|_{L^{\infty}}^2\|\nabla d\|_{H^3}^2+\|\nabla
u\|_{L^{\infty}}\|\nabla d\|_{H^3}^2,
\end{split}
\eeq
where $\epsilon$ will be chosen below. Since $|\nabla
d|^2+d\cdot\tr d=0$, we have
\beq\label{blp2.29}\begin{split}
J_4&=\int_{\R^3}\left( -D^{s}(|\nabla d|^2) d\cdot D^{s+2}d-D^{s}(|\nabla d|^2) D d\cdot D^{s+1}d \right)\,dx\\
&+\int_{\R^3}\left(D^{s}(|\nabla d|^2) D d\cdot D^{s+1}d+D^{s-1}(|\nabla d|^2) D^2 d\cdot D^{s+1}d \right)\,dx\\
&+\int_{\R^3}\left(D(|\nabla d|^2) D^{s} d\cdot D^{s+1}d+|\nabla d|^2 |D^{s+1}d|^2 \right)\,dx\\
&\leq C\int_{\R^3}\left(|\nabla^{4}d||\nabla d|+|\nabla^{2}d||\nabla^{3}d|\right)|D^{s+2}d|\,dx\\
&+ C\int_{\R^3}\left(|\nabla^{2}d|^3|\nabla^4 d|+|\nabla d||\nabla^{2}d||\nabla^{3}d||\nabla^{4}d|+|\nabla d|^2 |D^{s+1}d|^2 \right)\,dx\\
&\leq
\frac{1}{2}\|D^{s+2}d\|_{L^2}^2+C\int_{\R^3}(|\nabla^{4}d|^2|\nabla
d|^2+|\nabla^{2}d|^2|\nabla^3 d|^2)\,dx
+\|\nabla^2 d\|_{L^{6}}^3\|\nabla^4 d\|_{L^2}\\
&\leq \frac{1}{2}\|D^{s+2}d\|_{L^2}^2+C\|\nabla
d\|_{L^{\infty}}^2\|\nabla d\|_{H^3}^2.
\end{split}
\eeq
Combining (\ref{blp2.27}), (\ref{blp2.28}), with (\ref{blp2.29}), we have
 \beq\label{blp2.30}
 \frac{d}{dt}\|D^{s+1} d\|_{L^2}^2+\|D^{s+2} d\|_{L^2}^2
\leq\epsilon\|u\|_{H^4}^2 + C\left(\|\nabla
u\|_{L^{\infty}}+\|\nabla
d\|_{L^{\infty}}^2\right)\left(\|u\|_{H^3}^2+\|\nabla
d\|_{H^3}^2\right).
\eeq
Combining (\ref{blp2.26}) with (\ref{blp2.30}), we have
\bex\begin{split}
&\frac{d}{dt}(\|D^{s+1} d\|_{L^2}^2+\|D^{s}u\|_{L^2}^2)+\|D^{s+2} d\|_{L^2}^2+\|D^{s+1}u\|_{L^2}^2\\
\leq& \epsilon\|u\|_{H^4}^2+C\left(\|\nabla
u\|_{L^{\infty}}+\|\nabla
d\|_{L^{\infty}}^2\right)\left(\|u\|_{H^3}^2+\|\nabla
d\|_{H^3}^2\right).
\end{split}\eex
We can prove similar inequalities for all $|s|<3$. Summing over all
$s$ with $|s|\leq3$, and taking $\epsilon$ small enough, we have
\beq\label{blp2.31}
\frac{d}{dt}(\|\nabla d\|_{H^3}^2+\|u\|_{H^3}^2)
\leq C\left(\|\nabla u\|_{L^{\infty}}+\|\nabla d\|_{L^{\infty}}^2\right)\left(\|u\|_{H^3}^2+\|\nabla d\|_{H^3}^2\right).
\eeq

We now end our argument as follows.
Set $$m(t)=e+\|u\|_{H^3}+\|\nabla d\|_{H^3}.$$
Then by (\ref{blp2.31}), we have
$$\frac{dm(t)}{dt}\leq C\left(\|\nabla u\|_{L^{\infty}}+\|\nabla d\|_{L^{\infty}}^2\right)m(t).$$
By Gronwall's inequality, we obtain,
\beq\label{blp3.1}
m(t)\leq m(0)\exp\left(C\int^t_0(\|\nabla u(t)\|_{L^{\infty}}+\|\nabla d(t)\|_{L^{\infty}}^2)dt\right).
\eeq
By combining the following critical Sobolev embedding inequality (see \cite{Beale-Kato-Majda} for the detail)
\bex
\|\nabla u\|_{L^{\infty}}\leq C\left(1+\|\o\|_{L^2}+\|\o\|_{L^{\infty}}\ln(e+\|u\|_{H^3})\right),
\eex
with (\ref{blp2.22}), we have
$$\|\nabla u\|_{L^{\infty}}
\leq C\left(1+\|\o\|_{L^{\infty}}\ln(e+\|u\|_{H^3})\right).$$
Combining this inequality with (\ref{blp3.1}) and the ineqaulity $\ln(m(t))\geq 1$, we have,
\bex\begin{split}
\frac{d}{dt}\ln(m(t))&\leq C(1+\|\nabla d\|_{L^{\infty}}^2)+C\|\o\|_{L^{\infty}}\ln(m(t))\\
&\leq C(1+\|\nabla d\|_{L^{\infty}}^2+\|\o\|_{L^{\infty}})\ln(m(t)).
\end{split}\eex
By Gronwall's equality, we have
$$\ln(m(t))\leq \ln(m(0))\exp\left(C\int_0^t(\|\nabla d(t)\|_{L^{\infty}}^2+\|\o(t)\|_{L^{\infty}})dt\right),$$
Or equivalently,
$$m(t)\leq \exp\left(\ln(m(0))\exp\left(C\int_0^t(\|\nabla d(t)\|_{L^{\infty}}^2+\|\o(t)\|_{L^{\infty}})dt\right)\right),$$
for any $0\leq t\leq T_*$. This completes the proof.\endpf

\vspace{5mm}
{\bf Proof of Corollary \ref{2-d case}}: Assume that (\ref{blpcondition2-d}) were not true. Then
there is  $0<M_1<\infty$ such that
\beq\label{blp3.2}
\int_{0}^{T_*}\|\nabla d(t)\|^2_{L^{\infty}(\R^2)} dt\leq M_1.
\eeq
Then lemma \ref{blplemma2.1} implies
\beq\label{blp3.2}
\sup\ls_{0\leq t\leq T_*}\left(\|\o(t)\|_{L^{2}(\R^2)}+\|\nabla^2 d(t)\|_{L^{2}(\R^2)} \right)\leq C_1,
\eeq
where $C_1>0$ depends only on $u_0$, $d_0$ and $M_1$.
In particular, we have
$$(\nabla u,\nabla^2d)\in L^2_tL^2_x(\mathbb R^2\times [0,T_*]).$$
On the other hand, since $(u,d)$ satisfies (1.5), the following energy inequality holds (cf. \cite{Lin-Lin-Wang}):
\begin{equation}
\begin{split}
&\int_{\mathbb R^2}(|u(t)|^2+|\nabla d(t)|^2)\,dx
+2\int_0^t\int_{\mathbb R^2}(|\nabla u|^2+|\Delta d+|\nabla d|^2 d|^2)\,dxdt\\
&=\int_{\mathbb R^2}(|u_0|^2+|\nabla d_0|^2)\,dx
\end{split}
\end{equation}
for any $0<t\le T_*$. Therefore, we have
$$(u,\nabla d)\in L^\infty_tL^2_x(\mathbb R^2\times [0,T_*])\cap L^2_t H^1_x(\mathbb R^2\times [0,T_*]).$$
Applying the regularity theorem 1.2 of \cite{Lin-Lin-Wang}, we conclude $(u,d)\in C^\infty(\mathbb R^2\times (0,T_*])$.
This contradicts the assumption that $0<T_*<\infty$ is the first singular time. The proof is complete. \endpf

\end{document}